\documentclass{article}

\usepackage{arxiv}

\usepackage[utf8]{inputenc} 
\usepackage[T1]{fontenc}    
\usepackage{hyperref}       
\usepackage{url}            
\usepackage{booktabs}       
\usepackage{amsfonts}       
\usepackage{nicefrac}       
\usepackage{microtype}      
\usepackage{lipsum}
\usepackage{amsthm}
\usepackage[fleqn]{amsmath}
\usepackage[ruled]{algorithm2e}
\usepackage{multirow}
\usepackage{subfigure}
\usepackage{epsfig}

\title{Exact Approaches for Competitive Facility Location with Discrete Attractiveness}

\author{
  Yun Hui Lin \\
  Department of Industrial Systems Engineering \& Management\\
  National University of Singapore\\
  \texttt{isemlyh@gmail.com} \\
   \And
 Qingyun Tian\\
  School of Civil and Environmental Engineering\\
  Nanyang Technological University\\
  \texttt{qytian@ntu.edu.sg} \\
}

\begin{document}
\maketitle

\begin{abstract}
We study a variant of the competitive facility location problem, in which a company is to locate new facilities in a market where competitor's facilities already exist.  We consider the scenario where only a limited number of possible attractiveness levels is available, and the company has to select exactly one level for each open facility.  The goal is to decide the facilities' locations and attractiveness levels that maximize the profit. We apply the gravity-based rule to model the behavior of the customers and formulate a multi-ratio linear fractional 0-1 program. Our main contributions are the exact solution approaches for the problem. These approaches allow for easy implementations without the need for designing complicated algorithms and are “friendly” to the users without a solid mathematical background. We conduct computational experiments on the randomly generated datasets to assess their computational performance. The results suggest that the mixed-integer quadratic conic approach outperforms the others in terms of computational time. Besides that, it is also the most straightforward one that only requires the users to be familiar with the general form of a conic quadratic inequality. Therefore, we recommend it as the primary choice for such a problem.
\end{abstract}

\keywords{Competitive Facility Location \and Gravity Model \and Conic Programming \and  Outer Approximation \and Mixed-integer Linear Programming}

\section{Introduction}
The classical competitive facility location problem (CFL) studies a ``newcomer" company who enters a market where competitor's facilities already exist. The company has to decide the locations of its new facilities to maximize the market share or profit. In such a problem, customer's behaviors are typically considered. It is assumed customers will try to maximize their utility when selecting which facility to visit and their behaviors are interpreted probabilistically due to the randomness of customer's unobserved characteristics. Traditionally, the gravity model (also referred to as the Huff model) is used to forecast the probabilistic choice~\cite{drezner2002solving,drezner2018competitive,fernandez2007solving,kuccukaydin2011discrete}. Based on the model, the probability that a customer patronizes a facility is proportional to the attractiveness of the facility and declines according to a distance decay function $g(d)$. Here, the decay function represents the decline in the probability that a customer patronizes a facility as a function of his/her distance from it. A common-used  decay function takes the form of $g(d) = 1/d^{2}$. Alternatively, the multinomial logit model (MNL) has been applied to various CFLs~\cite{benati2002maximum,haase2014comparison,ljubic2018outer}.

There are two research streams on the CFL. The first one assumes that the competitor's facility location does not change when the newcomer makes decisions~\cite{benati2002maximum,drezner2002solving,fernandez2007solving,haase2014comparison,kuccukaydin2011discrete,lanvcinskas2017improving,ljubic2018outer}, while the other takes the perspective of the leader-follower game that considers the reaction of the competitor and studies the problem as a bi-level model ~\cite{drezner2015leader,kuccukaydin2011competitive,kuccukaydin2012leader,saiz2009branch}. We focus on the first stream. Specifically, we consider a variant of the classical single-level CFL, where additional decisions on the facility's attractiveness have to be made. By assuming that the attractiveness variables are continuous, related works (e.g.~\cite{fernandez2007solving,kuccukaydin2011discrete}) typically allow for unlimited variations in the possible value of attractiveness level. In this paper, we consider the case where only a limited number of possible attractiveness levels is available, and the company has to select exactly one for each open facility. This corresponds to the actual situation where only discrete-valued design options are present. We call this problem \textit{Competitive Facility Location with Discrete Attractiveness} (CFLD). Due to the integrity of attractiveness, this problem is more challenging to solve than the problems with continuous attractiveness.

To date, many approaches have been proposed to solve the CLF problem under the gravity-based rule. Most of them are based on metaheuristics or Lagrangian heuristic~\cite{drezner2002solving,kuccukaydin2011discrete,kuccukaydin2012leader,qi2017competitive}. For example, \cite{drezner2002solving} proposed five heuristic procedures and concluded that a two-step procedure that combines simulated annealing and an ascent algorithm yields the best solutions in general. \cite{kuccukaydin2011discrete} proposed three methods, namely, a Lagrangian heuristic, a branch-and-bound method with Lagrangian relaxation, and another branch-and-bound method with nonlinear programming relaxation. Their computational results show that the last method outperforms the others in terms of accuracy and efficiency. However, for the above approaches (not necessarily exact) to work well, one needs to design problem-dependent (even parameter-dependent) metaheuristics or implement a specialized branch-and-bound algorithm that makes use of Lagrangian relaxation bounds. This requires extensive coding efforts for computational implementation and is hard to generalize to more general cases. In this paper, we aim to provide straightforward and exact solution approaches that allow solving the problem without the need for designing complicated algorithms.

The contributions of this paper are that we first propose three easy-to-implement approaches for (exactly) solving CFLD and then conduct numerical experiments to investigate their performance on various problem instances. Our first approach is to reformulate the model as a mixed-integer linear program (MILP). However, the MILP approach requires a large number of new variables and constraints. The numerical experiment suggests that it cannot work well even for small-scale problems. The second approach is the outer approximation algorithm (OA) where we solve a sequence of (sparse) MILPs to find the exact solution. It has been applied to a similar problem where the attractiveness is not a decision variable~\cite{ljubic2018outer}. We show that it is also applicable to our problem. Finally, we show that the CFLD can be easily cast as a mixed-integer conic quadratic program (MICQP) and directly solved by second-order conic solvers. Recently, we have witnessed significant research contributions to (mixed) integer conic optimization~\cite{atamturk2011lifting,ccezik2005cuts,vielma2017extended}. The efficiency of the commercial solvers for MICQP has been undergoing rapid improvement as well. For example, Gurobi claims 2.8x improvement from version 6.0 to 6.5 \footnote{\url{https://www.gurobi.com/pdfs/WebSeminar-Gurobi-6.5.pdf}}. These developments make MICQP a promising approach for various problems. However, to the best of our knowledge, this approach has not been used to address the CFL. We demonstrate that this novel idea is straightforward to implement and exhibits computational advantages over MILP and OA. 

The rest of the paper is organized as follows. We introduce the model in Section~\ref{S:pd}. Section~\ref{S:solution} covers three solution approaches. We conduct numerical studies in Section~\ref{S:numerical} and conclude the paper in Section~\ref{S:conclusion}.

\section{Problem Description}
\label{S:pd}

\begin{table}
\caption{Nomenclature}
\label{tab:Nomenclature}       
\centering
\begin{tabular}{lll}
\hline\noalign{\smallskip}
\multicolumn{3}{l}{\underline{Sets}}\\
$I$ &      $:$ & set of customer zones \\
$J$ &     $:$ & set of candidate facilities\\
$K$ &     $:$ & set of competitor's facilities\\
$R$ &     $:$ & set of attractiveness level\\

\multicolumn{3}{l}{\underline{Parameters}}\\
$a_i$ & $:$ & anual buying power at customer zone $i$, $\forall i \in I$\\
$v_i$ &     $:$ & total utility of competitor's facilities for customer zone $i$, $\forall i \in I$\\
$d_{im}$ &     $:$ & distance between facility $m$ and customer zone $i$, $\forall i \in I, m \in  J \cup K$\\
$q_k$ & $:$ & attractiveness of competitor's facility $k$, $\forall k \in K$\\
$Q_{r}$ & $:$ & attractiveness of a facility at level $r$, $\forall r \in R$  \\
$f_j$ & $:$ & fixed cost of opening facility $j$, $\forall j \in J$\\
$c_{jr}$ & $:$ & attractiveness cost of facility $j$ at level $r$, $\forall j \in J, r \in R$\\

\multicolumn{3}{l}{\underline{Decision Variables}}\\
$x_j$ & $:$ &  boolean. If facility $j$ is open; 0, otherwise, $\forall j \in J$ \\
$y_{jr}$ & $:$ & boolean. If facility $j$ is open at level $r$; 0, otherwise, $\forall j \in J, r \in R$ \\
\noalign{\smallskip}\hline
\end{tabular}
\end{table}

Table~\ref{tab:Nomenclature} summarizes the main notations. In this paper, we use the gravity model to forecast customer's behaviors. We consider a company that wants to maximize its profit by determining the optimal location and attractiveness level of the facilities. We assume that the demands are aggregated at certain zones, denoted by set $I$, and the candidate facilities of the company are located at predetermined sites, denoted by set $J$. Meanwhile, there are existing facilities, denoted by set $K$, that belong to a competitor in the market. Using the Huff's gravity-based rule, the total utility of competitor's facilities for customers at zone $i$ is given by $v_i = \sum_{k \in K} q_k/d_{ik}^{2}, \forall i \in I$, where $d_{ik}$ is the distance between demand zone $i$ and existing facility $k$.
 
For the company,  when a facility is open, there are $|R|$ possible attractiveness levels that can be chosen with different costs. The attractiveness of facility $j$ is $A_j = \sum_{r \in R}Q_{r}y_{jr}, \forall j \in J$. According to the gravity-based rule,  the probability that customers at zone $i$ patronize facility $j$ is
\begin{equation}
p_{ij} =  \frac{A_j / d_{i j}^{2}}{\sum_{j \in J}\left( A_j / d_{i j}^{2}\right)+v_{i}} = \frac{ \sum_{r}Q_{r}y_{jr} / d_{i j}^{2}}{\sum_{j \in J}\left( \sum_{r}Q_{r}y_{jr} / d_{i j}^{2}\right)+v_{i}}
\end{equation}
Therefore, we compute the total revenue of the new facilities by
\begin{equation}
\sum_{i \in I} \sum_{j \in J} a_i p_{ij} =  \sum_{i \in I}^{n} a_{i} \frac{\sum_{j \in J}\left( \sum_{r}Q_{r}y_{jr} / d_{i j}^{2}\right)}{\sum_{j \in J}\left( \sum_{r}Q_{r}y_{jr} / d_{i j}^{2}\right)+v_{i}}
\end{equation}
To simplify the notation, we define $b_{ijr} = Q_{r}/(d_{ij}^{2} \cdot v_i)$. We can then formulate CFLD as the following integer nonlinear program:
\begin{equation} \label{p1}
\begin{aligned} 
\max~& \varphi(x, y)=\sum_{i \in I} a_{i} \frac{\sum_{j \in J} \sum_{r \in R} b_{i j r} y_{j r}}{\sum_{j \in J} \sum_{r \in R} b_{i j r} y_{j r}+1}-\sum_{j \in J} f_{j} x_{j}-\sum_{j \in J} \sum_{r \in R} c_{j r} y_{j r} \\ 
\text {st.} & \sum_{r \in R} y_{j r}=x_{j}, \forall j \in J \\ 
& x_{j} \in\{0,1\}, \forall j \in J\\
& y_{j r} \in\{0,1\}, \forall j \in J, r \in R 
\end{aligned}
\end{equation}
where the objective is to maximize the expected profit. The first term is the expected revenue collected by the new facilities. The second and third term are the fixed cost of opening the new facilities and the cost of operating the facilities in their attractiveness levels. The first constraint imposes that, if facility $j$ is not open, the attractiveness level of facility $j$ should be zero. On the other hand, if facility $j$ is open, there could be only one attractiveness level. The last two constraints are the binary restrictions on $x_j$ and $y_{jk}$.

\section{Solution Approach}
\label{S:solution}

In this section, we present three exact approaches. We first define the set:
\begin{equation}
\Omega = \left\{(x,y) :  \sum_{r \in R} y_{j r}=x_{j}, x_j\in\{0,1\}, y_{jr} \in\{0,1\}, \forall j \in J, r\in R \right\}
\end{equation}
which is the feasible region of (\ref{p1}). To maximize $\varphi(x,y)$ is equivalent to minimize $\sum_{i \in I} a_i - \varphi(x,y)$. We can write (\ref{p1}) as
\begin{equation}
\label{eqt:p1_obj} \min_{(x,y) \in \Omega}~\sum_{j \in J} f_{j} x_{j}+\sum_{j \in J} \sum_{r \in R} c_{j r} y_{j r}+\sum_{i \in I} \frac{a_{i}}{\sum_{j \in J} \sum_{r \in R} b_{i j} y_{j r}+1}
\end{equation}
where the last term in (\ref{eqt:p1_obj}) is a multi-ratio linear fractional term.  Therefore, (\ref{eqt:p1_obj})  belongs to the class of multi-ratio linear fractional 0-1 program (MLFP)~\cite{borrero2017fractional}.  In effect, (\ref{eqt:p1_obj}) is a special case of MLFPs because its continuous relaxation is a convex optimization problem. To support this argument, we define function $F_i(y)$ such that
\begin{align}
F_i(y) =  \frac{\sum_{j \in J}\sum_{r \in R}b_{ijr} y_{jr}}{\sum_{j \in J}\sum_{r \in R}b_{ijr} y_{jr}+ 1}, \forall i \in I
\end{align}
which is a concave function because it is the composition of a concave and increasing function, $z/(z+1)$ for $z \geq 0$, with an affine mapping $\sum_{j \in J}\sum_{r \in R}b_{ijr} y_{jr}$ \cite{boyd2004convex}. Therefore, the objective function of (\ref{eqt:p1_obj}) is convex. (\ref{eqt:p1_obj}) can then be studied as a mixed integer convex program and solved by convex optimization toolboxes. A popular one is CVXPY, which is a Python-embedded modeling framework for (mixed integer) convex optimization problems~\cite{cvxpy}. CVXPY automatically canonicalizes disciplined convex programs (DCP)~\cite{grant2006disciplined} to cone programs by expanding every nonlinear function into its graph implementation~\cite{grant2008graph}. Using CVXPY, we can call Gurobi, an advanced solver, to solve the DCP.  

Even with the convexity of the objective function, (\ref{eqt:p1_obj}) is still a hard problem due to the integrality constraint and the nonlinearity. In light of this, we present 3 alternative exact solution approaches.

\subsection{Mixed-integer Linear Program}

The MILP reformulation approach is widely used for MLFP~\cite{haase2014comparison,mendez2014branch} due to its ease of application and the powerfulness of the modern MILP solvers. Let $\beta_i = 1/(\sum_{j \in J}\sum_{r \in R}b_{ijr} y_{jr}+ 1), \forall i \in I$, then (\ref{eqt:p1_obj}) is equivalent to
\begin{equation}
\begin{aligned} \min~& \sum_{j \in J} f_{j} x_{j}+\sum_{j \in J} \sum_{r \in R} c_{j r} y_{j r}+\sum_{i \in I} a_{i} \beta_{i} \\ \text {st.}~& 1=\sum_{j \in J} \sum_{r \in R} b_{i j r} y_{j r} \beta_{i}+\beta_{i}, \forall i \in I \\ 
&(x,y) \in \Omega \end{aligned}
\end{equation}
which can be reformulated as a MILP by linearizing the bilinear term $ y_{jr} \beta_i$. Define $z_{ijr} = y_{jr} \beta_i$.  According to~\cite{adjiman1998global}, the convex lower bounds on $z_{ijr}$ are
\begin{equation}\label{constr:MC-1}
\begin{array}{l}{z_{i j r} \geq \beta_{i}-\beta_{i}^{U}\left(1-y_{j r}\right), \forall i \in I, j \in J, r \in R} \\ {z_{i j r} \geq \beta_{i}^{L} y_{j r}, \forall i \in I, j \in J, r \in R}\end{array}
\end{equation}
where $\beta^{U}_i$ is the upper bound of $\beta_i$. It can be set to 1. $\beta^{L}_i$ is the lower bound of $\beta_i$. It can be set to $\beta^{L}_i = 1/(\sum_{j \in J} b_{ij*} + 1)$, where $b_{ij*} = \max_{r \in R} b_{ijr}, \forall i \in I, j \in J$.  In addition, upper bounds on $z_{ijr}$ can be imposed to construct a better reformulation~\cite{mccormick1976}. This is achieved by adding two linear constraints:
\begin{equation}\label{constr:MC-2}
\begin{array}{l}{z_{i j r} \leq \beta_{i}-\beta_{i}^{L}\left(1-y_{j r}\right), \forall i \in I, j \in J, r \in R} \\ {z_{i j r} \leq \beta_{i}^{U} y_{j r}, \forall i \in I, j \in J, r \in R}\end{array}
\end{equation}
It is easy to see that (\ref{constr:MC-1})-(\ref{constr:MC-2}) enforce $z_{ijr} = y_{jr} \beta_i$. Therefore, (\ref{eqt:p1_obj}) can be expressed as the following MILP:
\begin{equation}
\begin{aligned}\min~ & \sum_{j \in J} f_{j} x_{j}+\sum_{j \in J} \sum_{r \in R} c_{j r} y_{j r}+\sum_{i \in I} a_{i} \beta_{i} \\ \text { st. } 1 &=\sum_{j \in J} \sum_{r \in R} b_{i j r} z_{i j r}+\beta_{i}, \forall i \in I \\ 
&(\ref{constr:MC-1})-(\ref{constr:MC-2}) \\
& (x,y) \in \Omega \end{aligned}
\end{equation}
which can be directly handled by MILP solvers. However, we need to introduce $O(|I|\cdot|J|\cdot|R|)$ variables and constraints. For large-scale problems, the resulting MILP may be difficult to solve to global optimality in a reasonable computational time.

\subsection{Outer Approximation}

The outer approximation (OA) has been applied to convex (concave) MLFPs~\cite{ljubic2018outer}. It proceeds by solving a sequence of sparse MILPs with an increasing number of constraints that refine the feasible region. The OA guarantees the convergence to global optimality in a finite number of iterations and it is efficient for large-scale MLFPs~\cite{elhedhli2005exact,ljubic2018outer}. 

Now, let $\hat{F}_i(y) = 1-F_{i}(y)$. We write (\ref{eqt:p1_obj}) in the epigraph form of $\hat{F}_i(y)$:
\begin{equation} \label{hygrah-p}
\begin{aligned} 
\min~& \sum_{j \in J} f_{j} x_{j}+\sum_{j \in J} \sum_{r \in R} c_{j r} y_{j r}+\sum_{i \in I} a_{i} \beta_{i} \\ 
\text { st. } & \beta_{i} \geq \hat{F}_{i}(y), \forall i \in I \\ 
&(x,y) \in \Omega \end{aligned}
\end{equation}
Given any point $\bar{y}$, since $\hat{F}_{i}(y)$ is a convex function, we can construct its lower bound by the first-order linear approximation on $\bar{y}$. The linear function will not eliminate any feasible region of (\ref{hygrah-p}). The following constraint is thus valid:
\begin{align}
\label{eqt:valid}  \beta_i \geq  \sum_{j \in J} \sum_{r \in R} \frac{\partial \hat{F}_{i}(\bar{y})}{\partial y_{jr}}(y_{jr} - \bar{y}_{jr}) + \hat{F}_{i}(\bar{y}),\forall i \in I
\end{align}
where $\frac{\partial \hat{F}_{i}(\bar{y})}{\partial y_{jk}}$ is the partial derivative of $\hat{F}_{i}$ with respect to $y_{jk}$. It is evaluated at $\bar{y}$. We can then model (\ref{hygrah-p}) using the following MILP formulation with $|I|$ additional variables and an exponential number of constraints:
\begin{equation} \label{OA:MILP}
\begin{aligned} \min~& \sum_{j \in J} f_{j} x_{j}+\sum_{j \in J} \sum_{r \in R} c_{j r} y_{j r}+\sum_{i \in I} a_{i} \beta_{i} 
\\ \text {st.}~& \beta_{i} \geq \sum_{j \in J} \sum_{r \in R} \frac{\partial \hat{F}_{i}(\bar{y})}{\partial y_{j r}}\left(y_{j r}-\bar{y}_{j r}\right)+\hat{F}_{i}(\bar{y}), \forall i \in I, \bar{y} \in T
\\ &(x,y) \in \Omega \end{aligned}
\end{equation}
The details of the OA is presented in Algorithm~\ref{alg:OA}. Here, $T$ is the set of recorded points $(x,y)$. Each point (expect the initial one) is obtained by solving the MILP during each iteration. It defines a valid inequality that is further added to the MILP. As the number of points increases, (\ref{OA:MILP}) provides better approximation for (\ref{hygrah-p}) and ultimately converges to the optimal solution (when it finds a solution that already exists in set $T$) in a finite number of iterations. 

\begin{algorithm}[H] \label{alg:OA}
\SetAlgoLined
Step 0: Initialize $x^0$ and $y^0$; $n=0$; $T = \emptyset$.\\
Step 1: $T: = T \cup (x^n,y^n)$. Compute the partial derivative:
$$\frac{\partial \hat{F}^n_i}{\partial y_{jr}} :=  \frac{-b_{ijr}}{(\sum_{j \in J}\sum_{r \in R}b_{ijr} y^n_{jr}+ 1)^2}, \forall i \in I, j \in J, r \in R$$\\
Step 2: Solve (\ref{OA:MILP}) to obtain the solution $(x^{n+1},y^{n+1})$. \\
Step 3: If $(x^{n+1},y^{n+1}) \in T$, stop and output $(x^{n+1},y^{n+1})$. Else,  $n := n +1$ and go to Step 1.
 \caption{Outer Approximation Algorithm}
\end{algorithm}

\subsection{Mixed-integer Conic Quadratic Program}
\label{S:micqp}
Next, we show that (\ref{eqt:p1_obj}) can be recast as a mixed-integer conic quadratic program (MICQP) that can be efficiently solved by off-the-shelf solvers (e.g. CPLEX, Gurobi, Mosek). Conic optimization refers to optimizing a linear function over conic inequalities~\cite{ben2001lectures}. A conic quadratic inequality take the form of
\begin{equation}\label{conic_form}
|| Ax - b||_2 \leq c^Tx - d
\end{equation}
where $||\cdot||_2$ is the L2 norm. $A$, $b$, $c$ and $d$ are matrices or vectors with proper sizes. Generally, a rotated cone on $x,y,z \geq 0$:
\begin{align}
x^2 \leq yz
\end{align}
can be represented as an equivalent conic quadratic inequality:
\begin{align}
|| 2x, y-z||_2 \leq y + z
\end{align}
We make use of the rotated cone in our reformulation.

Let $z_i = \sum_{j \in J}\sum_{r \in R}b_{ijr} y_{jr}+ 1$, then $\hat{F}_i(y) = 1/z_i$ and $z_i \geq 1$. Now, define $\beta_i$ such that $\beta_i \geq \hat{F}_i(y)$. We have the rotated cone inequality: $\beta_i z_i \geq 1$. Therefore, (\ref{eqt:p1_obj}) is equivalent to the following MICQP:
\begin{equation} \label{CQMIP}
\begin{aligned} \min~& \sum_{j \in J} f_{j} x_{j}+\sum_{j \in J} \sum_{r \in R} c_{j r} y_{j r}+\sum_{i \in I} a_{i} \beta_{i} \\ \text {st.}~& z_{i}=\sum_{j \in J} \sum_{r \in R} b_{i j r} y_{j r}+1, \forall i \in I \\ & \beta_{i} z_{i} \geq 1, \forall i \in I \\ & z_{i} \geq 1, \forall i \in I \\ & 0 \leq \beta_{i} \leq 1, \forall i \in I \\ &(x,y) \in \Omega \end{aligned}
\end{equation}
where the objective is linear and the constraints are either conic quadratic or linear. 
It can be verified that $\beta_i z_i = 1$ follows in the optimal solution. Therefore, we can enforce that $\beta_i$ is not larger than 1 since $z_i \geq 1$. To date, MICQP models have recently been employed to address various problems, such as portfolio optimization~\cite{vielma2008lifted}, location-inventory problem~\cite{atamturk2012conic}, and assortment optimization~\cite{csen2018conic}. However, to the best of our knowledge, this approach has not been used to address the CFL. As in (\ref{CQMIP}), this novel idea is straightforward to implement.  The MICQP contains only $|I|$ conic quadratic constraints and requires only $O(|I|)$ additional variables.

\section{Numerical Test}
\label{S:numerical}

To investigate the computational performance of three approaches, we use artificially generated data. We assume the number of demand zones and the number of candidate facilities are equal to each other, that is, $|I|=|J|$. The number of attractiveness levels $|R|$ is 5. Therefore, the number of binary variables is $6\cdot|I|$. Meanwhile, we set $Q = [100, 300, 500, 700, 900]$ and $c_{j\cdot} = 2*Q,\forall j \in J$. The attractiveness of  a competitor's facility $q_k$ is generated from a uniform distribution $[100,1000]$. For the distance matrix, the locations of the demand zones, the candidate facilities, and the competitor's facilities are generated from a 2-dimensional uniform distribution $[0,100]^2$. The distance is then calculated by the Euclidean distance. Finally, the annual buying power $a_i$ is generated from a uniform distribution $[100,10000]$.  All computational experiments are done using Python on a 16 GB memory iOS computer with a 2.6 GHz Intel Core i7 processor. We use Gurobi 8.1.1 under default settings as the solver. 

We start by testing MILP, OA and MICQP on small-scale instances, namely, $I=20$ and $I=30$. We set the maximal computational time to 3600s. Table~\ref{tab:time1} reports the results. In all instances, the computational times by OA and MICQP are minimal (less than 5s). The MILP approach, however, requires significant computational loads. In particular, when there are 30 customers, this approach fails to solve all instances in 1 hour. In effect, the average MIP relative gap for those unsolved instances in Table~\ref{tab:time1} are significant, 
but the objective value yielded by the MILP approach is optimal expect for the instance with $I=30, f=0, K=1$. This indicates that the continuous relaxation of the MILP approach is rather weak, which, in many cases, presents the solver from closing the MIP gap even if the optimal solution has been found and leads to extensive branching, thus, a long computational time. Therefore, OA and MICQP are better solution approaches and the MILP works poorly for CFLD, even for small-scale instances.

\begin{table}
\centering
\caption{Computational results of MILP, OA and MICQP on small-scale instances. CPU is the computational time (s); N.A. indicates the problem cannot be solved in $3600s$; $\varphi$ is the objective, measured in thousands; F is the number of new facilities.}  \label{tab:time1} \centering \scriptsize \label{tab:time1}
\setlength{\tabcolsep}{1.5mm}{
\begin{tabular}{rrrrrrrrrrrrrrrrr} 
\hline
 \multicolumn{1}{c}{I}  &  \multicolumn{1}{c}{$f$}    &  \multicolumn{1}{c}{K}  & \multicolumn{3}{c}{MILP} &  & \multicolumn{3}{c}{OA~} &  & \multicolumn{3}{c}{MICQP~}  \\ 
\cline{4-6}\cline{8-10}\cline{12-14}
   &      &    &  \multicolumn{1}{c}{CPU}    &  \multicolumn{1}{c}{$\varphi$}    &  \multicolumn{1}{c}{F}    &  &  \multicolumn{1}{c}{CPU} &  \multicolumn{1}{c}{$\varphi$}    &  \multicolumn{1}{c}{F}      &  &  \multicolumn{1}{c}{CPU} &  \multicolumn{1}{c}{$\varphi$}    &  \multicolumn{1}{c}{F }         \\ 
\hline
20 & 0    & 1  &    N.A.    &    114.4    &  9      &  &  1.6   &   114.4     &   9       &  &   0.3  &     114.4   &    9          \\
   &      & 5  &    N.A.    &   97.1     &    11    &  &   1.1  &   97.1     &     11     &  &     0.2 &   97.1     &     11         \\
   &      & 10 &   N.A.     & 85.9    & 12        &  & 0.9    &    85.9      &  12        &  &    0.2 &      85.9   &  12            \\
   &      &    &        &        &        &  &     &        &          &  &     &        &              \\
   & 2500 & 1  & 2215.9 & 103.6 & 3      &  & 0.8 & 103.6 & 3        &  & 0.4 & 103.6 & 3            \\
   &      & 5  & N.A.   & 81.0  & 5      &  & 1.1 & 81.0  & 5        &  & 0.3 & 81.0  & 5            \\
   &      & 10 & 2548.8 & 67.1  & 5      &  & 0.6 & 67.1  & 5        &  & 0.1 & 67.1  & 5            \\
   &      &    &        &        &        &  &     &        &          &  &     &        &              \\
   & 5000 & 1  & 256.0  & 96.1  & 2      &  & 1.0 & 96.1 & 2        &  & 0.2 & 96.1  & 2            \\
   &      & 5  & 339.0  & 70.7  & 4      &  & 0.6 & 70.7  & 4        &  & 0.1 & 70.7  & 4            \\
   &      & 10 & 110.1  & 56.2  & 4      &  & 0.2 & 56.2  & 4        &  & 0.1 & 56.2  & 4            \\
   &      &    &        &        &        &  &     &        &          &  &     &        &              \\
30 & 0    & 1  &  N.A.      &  147.5      &  12      &  &  4.6   & 147.6       &  13        &  &  0.7   &   147.6     &        13      \\
   &      & 5  &   N.A.     &  130.3      &   14     &  &  4.7   &  130.3      &     14     &  &  0.4    &   130.3     &     14         \\
   &      & 10 &  N.A.      &  103.3      &   15     &  &  1.4   &  103.3      &   15       &  & 0.3    &     103.3   &         15     \\
   &      &    &        &        &        &  &     &        &          &  &     &        &              \\
   & 2500 & 1  & N.A.   & 134.3 & 4      &  & 2.0 & 134.3 & 4        &  & 0.2 & 134.3 & 4            \\
   &      & 5  & N.A.   &  110.6     &   7     &  & 3.0 & 110.6 & 7        &  & 0.5 & 110.6 & 7            \\
   &      & 10 & N.A.   & 78.8      &  8      &  & 1.4 & 78.8  & 8        &  & 0.4 & 78.8  & 8            \\
   &      &    &        &        &        &  &     &        &          &  &     &        &              \\
   & 5000 & 1  &   N.A.     &  125.9      &  3      &  & 1.3 & 125.9       &  3        &  & 0.4 & 125.9       &    3          \\
   &      & 5  &    N.A.    &  96.6      &   5     &  & 5.0 &  96.6      &  5        &  & 0.6 &  96.6      &     5         \\
   &      & 10 &   N.A.     & 61.9       & 5       &  & 1.4 &   61.9     &  5        &  & 0.6 & 61.9       &     5         \\
\hline
\end{tabular}}
\end{table}

Next, we test OA and MICQP on larger-scale instances, namely, $I=75$ and $I=100$. To show the advantages of our approaches, we also solve the model by CVXPY. Using this toolbox, we can call Gurobi to solve the problem. Therefore, CVXPY is an appropriate benchmark for our approaches. We set the maximal computational time to 10000s. 

Table~\ref{tab:OA_MICQP} reports the computational time. We first look at the instances with 75 customer zones. All approaches can solve them to global optimality.  In general, CVXPY has a longer computational time than OA and MICQP, meaning that our approaches indeed facilitate the problem solving. The most efficient approach is MICQP. In all tested instances, the computational time required for MICQP is less than 70s. We then look at the instances with 100 customer zones. We observe that CVXPY cannot work well on these instances and only solves $8/18$ instances in 10000s. Among the approaches, MICQP is the only one that solves all instances. Compared to OA, the computational advantage of MICQP is significant, especially when the fixed costs are 2500 and 5000. Indeed, when OA is applied to large-scale instances, it generally takes about $6$ to $10$ iterations to converge and needs to solve the same number of large-scale MILPs. It may experience a significant slow-down at some iterations because the number of constraints increases during the iteration.

The numerical experiments suggest that MICQP outperforms CVXPY and OA in terms of computational time. In effect, as in Section~\ref{S:micqp}, it is easy to reformulate the CFLD model into a MICQP. Therefore, we recommend MICQP as the approach for solving CFLD.

\begin{table}[]\caption{Computational time (s) by CVXPY, OA and MICQP. N.A. indicates the problem cannot be solved in $10000s$; Avg is the average time; and $t^*$ indicates the Avg is at least $t$.}  \label{tab:OA_MICQP} \centering \scriptsize 
\setlength{\tabcolsep}{1.5mm}{ 
\begin{tabular}{rrrrrrrrrrrrrrrrrrrr}
\hline
\multirow{2}{*}{I} & \multirow{2}{*}{K} & \multicolumn{3}{c}{$f=2500$} &  & \multicolumn{3}{c}{$f=5000$}     &  & \multicolumn{3}{c}{$f=10000$} \\ \cline{3-5} \cline{7-9} \cline{11-13} 
                   &                    &  \multicolumn{1}{c}{CVXPY}    &  \multicolumn{1}{c}{OA}      &  \multicolumn{1}{c}{MICQP}   &  &  \multicolumn{1}{c}{CVXPY}  &  \multicolumn{1}{c}{OA}     &  \multicolumn{1}{c}{MICQP}          &  &  \multicolumn{1}{c}{CVXPY}      &  \multicolumn{1}{c}{OA}      &  \multicolumn{1}{c}{MICQP}  \\ \hline
75                 & 5                  & 2309.6   & 432.6   & 69.3    &  & 238.1  & 72.3   & 9.3            &  & 129.4      & 30.3    & 16.2   \\
                   & 10                 & 199.5    & 124.5   & 9.9     &  & 86.7   & 32.8   & 5.9            &  & 51.3       & 20.7    & 6.4    \\
                   & 15                 & 69.2     & 130.5   & 5.7     &  & 78.2   & 31.1   & 4.9            &  & 30.7       & 11.7    & 5.2    \\
                   & 20                 & 29.2     & 155.0   & 5.0     &  & 60.4   & 24.6   & 4.3            &  & 25.2       & 10.5    & 5.2    \\
                   & 25                 & 14.3     & 106.3   & 2.5     &  & 35.9   & 35.4   & 5.2            &  & 26.7       & 12.1    & 4.1    \\
                   & 30                 & 15.1     & 137.6   & 2.3     &  & 27.7   & 92.2   & 5.9            &  & 13.4       & 12.0    & 2.4    \\
                   & \textbf{Avg}       & 439.5    & 181.1   & 15.8    &  & 87.8   & 48.1   & 5.9            &  & 46.1       & 16.2    & 6.6    \\
                   &                    &          &         &         &  &        &        &                &  &            &         &        \\
100                & 10              & N.A.      & N.A.      & 3161.9  &  &  N.A.      & 1654.4 & 126.8           &  & N.A.         & 231.6   &168.1  \\
                   & 20                 & N.A.       & 5621.2  & 108.0   &  &   N.A.     & 1233.5 & 66.3          &  & N.A.         & 478.6   & 375.8  \\
                   & 30                 &  N.A.      & 596.7   & 17.1    &  &    N.A.    & 6178.8 & 516.7          &  & 5143.1     & 635.5   & 93.5  \\
                   & 40                 & 299.9    & 126.1   & 8.1     &  & N.A.       & 520.0  & 50.4           &  & 956.8      & 50.5    & 10.9  \\
                   & 50                 & 84.6     & 539.5   & 5.6     &  & N.A.     & 152.2  &  39.6           &  & 449.0      & 22.0    & 5.9    \\
                   & 60                 & 55.3     & 430.5   & 3.9     &  & 5310.1 & 140.9  & 19.1          &  & 48.1       & 9.9     & 4.7    \\
                   & \textbf{Avg}  & $5073.7^*$         & $2885.7^*$    &  550.8         &  & $9218.3^*$       & 1646.6 & 136.5          &  & $4433.8^*$  & 238.0   & 109.8   \\ \hline
\end{tabular}}
\end{table}

\begin{figure*}[h]
\begin{center}
	\subfigure[]{
      \psfig{figure=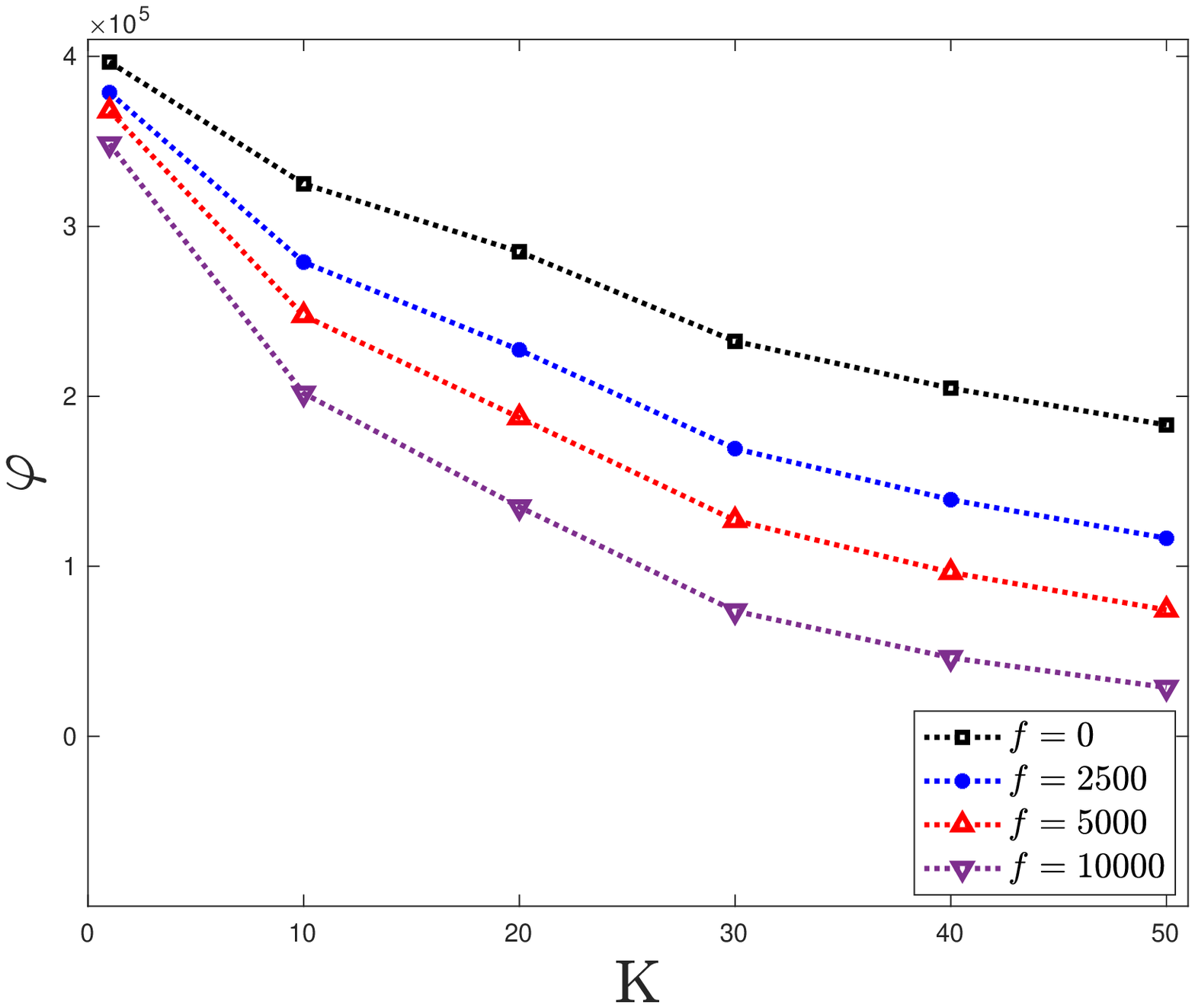,width=70mm}\label{fig:profit}
    }
    \subfigure[]{
      \psfig{figure=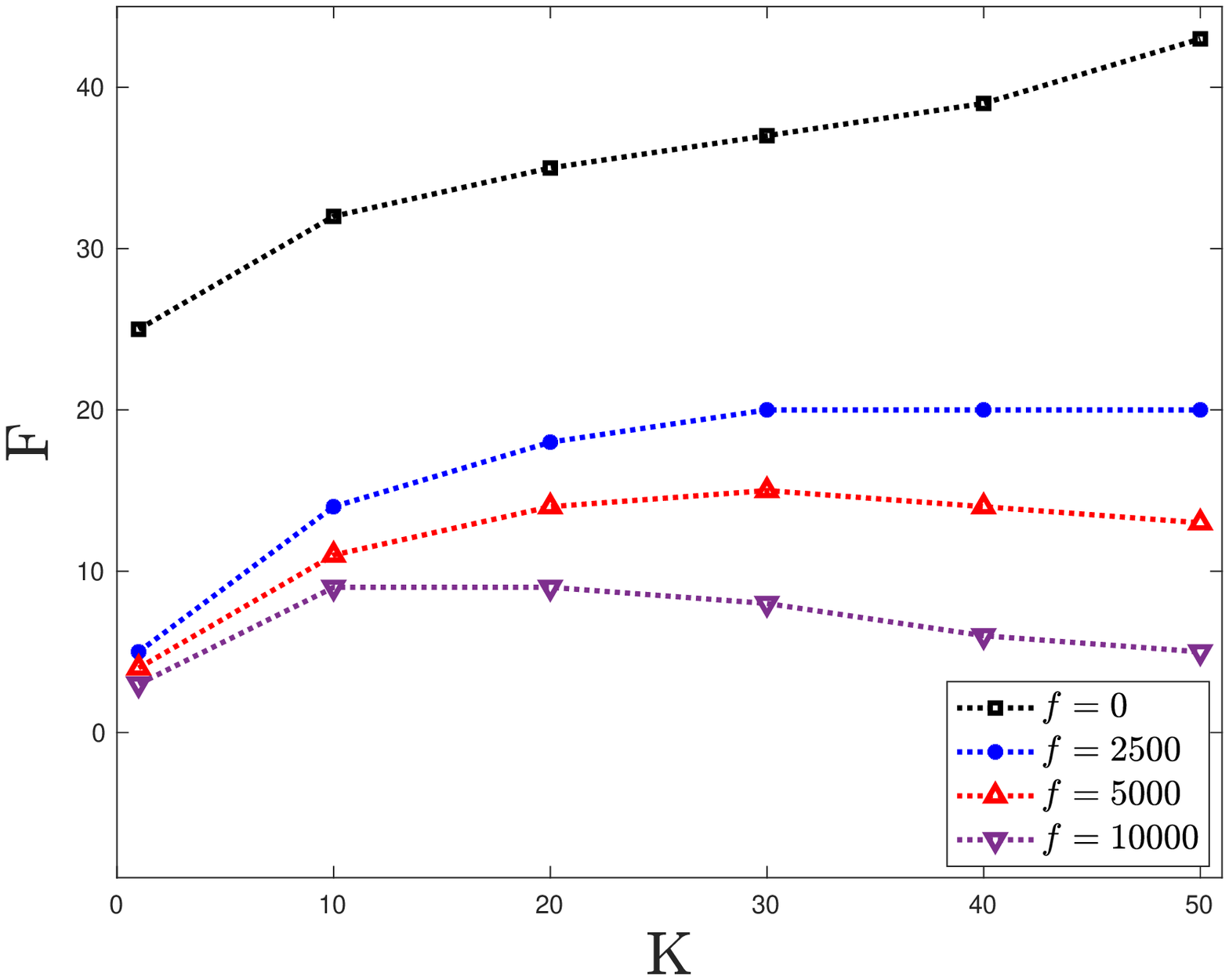,width=70mm}\label{fig:F}
      }  
\end{center}
\caption{Profit ($\varphi$) and the number of new facilities (F) under different parameter values.} \label{fig:K}
\end{figure*}

Finally, using the instances with $75$ customer zones, we briefly discuss the impacts of $K$ and $f$ on the solution and profit. The results are derived by MICQP. In Fig~\ref{fig:profit}, as the number of competitor's facilities (K) increases, the profit decreases. This is reasonable because the competitor will be able to attract more customers with more facilities. Meanwhile, the profit also intuitively decreases with $f$. Fig~\ref{fig:F} depicts the number of new facilities versus K under varying fixed costs. Given any K value, F is larger when $f$ is lower. It is within our expectation that with a lower fixed cost, the company will tend to open more facilities in order to attract more customers and increase its profit.  The interesting results in the figure come from the different trends of F versus the increasing K under different fixed costs. When $f$ is 0, F shows an obvious growing trend. The model suggests the company to aggressively open new facilities and compete for the demand with the competitor. When the fixed cost is minimal, the advantage of opening new facilities is more likely to overweight the additional cost. However, when $f$ is significant, if there exists a large number of facilities (i.e., a high level of competition), it is then better to remain conservative because the cost of opening additional facilities could be higher than the increased revenue. Therefore, when $f=10000$, F first increases with K and then decreases. 

\section{Conclusion}
\label{S:conclusion}

We studied the competitive facility location problem with a finite number of choices on the facility's attractiveness levels. This problem is inherently more challenging to solve than the problems with continuous attractiveness. We presented three exact approaches. These approaches are easy-to-implement and user-friendly. Our computational results highlighted that the MICQP approach is superior to others in terms of computational time. It also outperforms a popular convex optimization toolbox, CVXPY, in all of our tested instances. For the implementation, the MICQP approach only requires the users to be familiar with the general form of a conic quadratic inequality. Therefore, we recommend it as the primary choice for such a problem.  Given the promise of MICQP formulation for the CFL, it is worthwhile to explore its capabilities on other CFL problems based on different customer choice models (e.g. multinomial logit model) and to investigate its potential applications to the nonconvex CFL problems where OA and CVXPY are not applicable.

\bibliographystyle{plain}
\bibliography{references}
\end{document}